\documentclass[a4paper,12pt]{article}
\usepackage{amsfonts}
\newcommand{\grad}{\mathop{\rm grad}\nolimits}
\newcommand{\curl}{\mathop{\rm curl}\nolimits}
\newcommand{\intprod}{\;\rule{5pt}{.3pt}\rule{.3pt}{7pt}\;}
\renewcommand\thefootnote{$\ast$}
\begin{document}
\begin{center}
{\Large\bf Ricci curvature and the mechanics of solids}\\[16pt]
Michael Eastwood\footnote{Mathematical Sciences Institute, Australian National 
University, ACT 0200.}
\renewcommand\thefootnote{}
\footnotetext{Support from the Australian Research Council is gratefully 
acknowledged.}
\end{center}
\begin{abstract}
We discuss some differential geometry pertaining to continuum mechanics and the
route recently taken by D.N.~Arnold, R.S.~Falk, and R.~Winther in deriving new
improved finite element schemes in linear elasticity from constructions in
projective geometry.
\end{abstract}
\renewcommand\thefootnote{$\dagger$}

\subsection*{Some vector analysis}
We start with some basics from vector analysis~\cite{spain}. Let us write
$$\partial_1\equiv \partial/\partial x_1\qquad
\partial_2\equiv \partial/\partial x_2\qquad
\partial_3\equiv \partial/\partial x_3$$
for the partial derivatives in~${\mathbb{R}}^3$. The gradient of a smooth
function $f$ defined on $U^{\mathrm{open}}\subseteq{\mathbb{R}}^3$ is the 
vector field
$$\grad f
\equiv(\partial_1f,\partial_2f,\partial_3f)$$
on~$U$. If $X=(X_1,X_2,X_3)$ is a smooth vector field on~$U$, then
$$\curl X
\equiv(\partial_2X_3-\partial_3X_2,
\partial_3X_1-\partial_1X_3,\partial_1X_2-\partial_2X_1).$$
It is readily verified that $\curl\circ\grad=0$. Indeed, if 
$U$ is a sufficiently simple set, such as a ball, then
\begin{equation}\label{poincare}
\framebox{$X=\grad f,\mbox{ for some }f\iff\curl X=0$.}
\end{equation}

In effect, the curl of a vector field is the skew part of
the $3\times 3$ matrix $(\partial_iX_j)$ of partial derivatives.
Let us instead consider the symmetric part
$$\textstyle\Sigma=(\Sigma_{ij})\equiv(\frac12[\partial_iX_j+\partial_jX_i])$$
and ask for conditions that a given symmetric tensor field
$\Sigma=(\Sigma_{ij})$ be of this form. The answer is that $\Sigma$
should satisfy the {\em Saint-Venant\/}
equations $\curl\curl\Sigma=0$, where
$\curl\curl\Sigma$ is the symmetric matrix obtained by
\begin{itemize}
\item firstly regarding $\Sigma$ as a row vector (whose entries just happen to
be column vectors) to form $\curl\Sigma$, 
\item then regarding $\curl\Sigma$ as a column vector (whose entries just
happen to be row vectors) to form $\curl(\curl\Sigma)$.
\end{itemize}
The statement (\ref{poincare}) has a useful counterpart as
follows.
\begin{equation}\label{counterpart}
\framebox{$\Sigma_{ij}=
\frac12[\partial_iX_j+\partial_jX_i],\mbox{ for some }X
\iff\curl\curl\Sigma=0$.}
\end{equation}
Indeed, we shall see that (\ref{counterpart}) can be deduced from 
(\ref{poincare}). This, in turn, has consequences in the design of finite 
element schemes concerned with elasticity.

\subsection*{Ricci curvature in three dimensions}
Readers unfamiliar with differential geometry might omit this section on 
first reading.
There is also a close link between (\ref{counterpart}) and Ricci curvature in
three dimensions. 
Using the Einstein summation convention, if
$g_{ij}$ is a Riemannian metric with inverse~$g^{ij}$, then the
{\em Ricci\/} tensor $R_{ij}$ is the symmetric tensor given by 
$$\partial_k\Gamma_{ij}{}^k-\partial_i\Gamma_{jk}{}^k
+\Gamma_{ij}{}^m\Gamma_{mk}{}^k-\Gamma_{ik}{}^m\Gamma_{jm}{}^k,
\mbox{ where }
\textstyle\Gamma_{ij}{}^k\equiv
\frac12g^{kl}[\partial_ig_{jl}+\partial_jg_{il}-\partial_lg_{ij}].$$
In three dimensions $R_{ij}=0$ if and only if $g_{ij}$ is flat, meaning that
there is a local change of co\"ordinates that transforms $g_{ij}$ as a tensor
into the flat metric~$\delta_{ij}$ (more specifically, there is a change of
co\"ordinates with Jacobian matrix $J$ such that $(g_{ij})=J^tJ$). The
infinitesimal version of this statement is essentially (\ref{counterpart}).
More precisely, if $\Sigma_{ij}$ is an arbitrary symmetric tensor on
$U^{\mathrm{open}}\subseteq{\mathbb{R}}^3$ and we consider the metric
$g_{ij}^\epsilon=\delta_{ij}+\epsilon\Sigma_{ij}$ where $\epsilon$ is
sufficiently small that $g_{ij}^\epsilon$ is positive definite, then
$$(\curl\curl\Sigma)_{ij}=\frac{d}{d\epsilon}G^\epsilon_{ij}|_{\epsilon=0}$$ 
where $G_{ij}$ is the Einstein tensor $Rg_{ij}-2R_{ij}$ for $R=g^{kl}R_{kl}$.
The Einstein tensor carries the same information as the Ricci tensor but has
the advantage that the Bianchi identity simply says that $G_{ij}$ is
divergence-free $\nabla^iG_{ij}=0$. Sure enough, one can readily verify that 
$\partial^i(\curl\curl\Sigma)_{ij}=0$.

\subsection*{Translation from (\ref{poincare}) to (\ref{counterpart})}
Firstly, some convenient notation in three dimensions. Let us write
$\epsilon_{ijk}$ for the totally skew tensor with $\epsilon_{123}=1$. It 
allows us to write $(\curl X)_i=\epsilon_i{}^{jk}\partial_jX_k$. 
Now consider a pair $F=(X_\ell,Y_\ell)$ of vector fields on
$U^{\mathrm{open}}\subseteq{\mathbb{R}}^3$, regarded as a function with
values in the vector space
${\mathbb{W}}\equiv{\mathbb{R}}^3\oplus{\mathbb{R}}^3$. If we define the 
gradient of $F$ by
\begin{equation}\label{tractors}
\grad\left[\begin{array}cX_\ell\\ Y_\ell\end{array}\right]=
\left[\begin{array}c\partial_jX_\ell-\epsilon_{j\ell}{}^mY_m\\ 
\partial_jY_\ell\end{array}\right]\end{equation}
and use this definition na\"{\i}vely to compute the $\curl$ of a vector field 
with values in ${\mathbb{W}}$, then we obtain
$$\curl\left[\begin{array}c\Sigma_{j\ell}\\ \Xi_{j\ell}\end{array}\right]=
\left[\begin{array}c\epsilon_i{}^{jk}\partial_j\Sigma_{k\ell}
-\epsilon_i{}^{jk}\epsilon_{j\ell}{}^m\Xi_{km}\\ 
\epsilon_i{}^{jk}\partial_j\Xi_{k\ell}\end{array}\right]=
\left[\begin{array}c\epsilon_i{}^{jk}\partial_j\Sigma_{k\ell}
-\Xi_{\ell i}+\delta_{i\ell}\Xi_m{}^m\\ 
\epsilon_i{}^{jk}\partial_j\Xi_{k\ell}\end{array}\right].$$

It is readily verified that $\curl\circ\grad=0$. This says precisely that
(\ref{tractors}) defines a flat connection, which enables one to deduce that if
$U$ is a sufficiently simple set, such as a ball, then
\begin{equation}\label{coupledpoincare}
\framebox{$\Psi=\grad F,\mbox{ for some }F\iff\curl \Psi=0$.}
\end{equation}
To deduce~(\ref{counterpart}), let us suppose that $\Sigma_{ij}$ is symmetric 
and set
$$\Psi=\left[\begin{array}c\Sigma_{j\ell}\\ \Xi_{j\ell}\end{array}\right]=
\left[\begin{array}c\Sigma_{j\ell}\\ \epsilon_\ell{}^{im}\partial_i\Sigma_{mj}
\end{array}\right],\enskip\mbox{so that}\enskip\curl\Psi=
\left[\begin{array}c0\\ (\curl\curl\Sigma)_{i\ell}
\end{array}\right].$$
If $\curl\curl\Sigma=0$, we immediately infer the existence of vector fields
$X_\ell$ and $Y_\ell$ on $U$ such that
$$\grad\left[\begin{array}cX_\ell\\ Y_\ell\end{array}\right]=\Psi\quad
\mbox{i.e.}\quad
\left[\begin{array}c\partial_jX_\ell-\epsilon_{j\ell}{}^mY_m\\ 
\partial_jY_\ell\end{array}\right]=
\left[\begin{array}c\Sigma_{j\ell}\\ \epsilon_\ell{}^{im}\partial_i\Sigma_{mj}
\end{array}\right].$$
In particular, $\Sigma_{j\ell}=\frac12[\partial_jX_\ell+\partial_\ell X_j]$, as
required.

\subsection*{Continuum mechanics}
Although different words are used, Riemannian differential geometry in three
dimensions is exactly what is needed to set up the mechanics of
solids~\cite{ciarlet}. The metric tensor is known as the {\em strain\/} in
continuum mechanics. The Einstein tensor is known as the {\em stress}. The
Bianchi identity says that the stress tensor is divergence-free, interpreted as
a {\em conservation law\/} in mechanics. Linearising around the flat metric
gives the following complex of tensors on~${\mathbb{R}}^3$
\begin{equation}\label{linearelasticitycomplex}\begin{array}{ccccccc}
X_i&\mapsto&\frac12[\partial_jX_j+\partial_jX_i]&&S_{ij}&\mapsto&
\partial^iS_{ij}\\
\mbox{displacement}&\to&\mbox{strain}&\to&
\mbox{stress}&\to&\mbox{load}\,,\\
&&\Sigma_{ij}&\mapsto&\epsilon_i{}^{km}\epsilon_j{}^{\ell n}
\partial_k\partial_\ell\Sigma_{mn}
\end{array}\end{equation}
where the {\em displacement\/} and {\em load\/} are vector fields whilst the 
stress and strain are symmetric $2$-tensors.

\subsection*{Finite element schemes}
We have already seen that the flat connection (\ref{tractors}) somehow embodies
the operator $\Sigma\mapsto\curl\curl\Sigma$ relating strain and stress
in~(\ref{linearelasticitycomplex}). More generally and precisely, the whole
complex (\ref{linearelasticitycomplex}) may be derived from the
connection~(\ref{tractors}). To do this, recall that the gradient operator
(\ref{tractors}) concerned functions with values in
${\mathbb{W}}={\mathbb{R}}^3\oplus{\mathbb{R}}^3$. Thus, we may write
\begin{equation}\label{bigdiagram}\begin{array}{ccccccc}{\mathbb{W}}
&\stackrel{{\rm grad}}{\longrightarrow}&{\mathbb{R}}^3\otimes{\mathbb{W}}
&\stackrel{{\rm curl}}{\longrightarrow}&{\mathbb{R}}^3\otimes{\mathbb{W}}
&\stackrel{{\rm div}}{\longrightarrow}&{\mathbb{W}}\\
\|&&\|&&\|&&\|\\
{\mathbb{R}}^3&&
\framebox{${\mathbb{R}}^3$}\oplus{\mathrm{S}}^2{\mathbb{R}}^3
&&\framebox{${\mathbb{R}}^3\otimes{\mathbb{R}}^3$}&&
\framebox{${\mathbb{R}}^3$}\\
\oplus&\nearrow&\oplus&\nearrow&\oplus&\nearrow&\oplus\\
\framebox{${\mathbb{R}}^3$}
&&\framebox{${\mathbb{R}}^3\otimes{\mathbb{R}}^3$}
&&{\mathrm{S}}^2{\mathbb{R}}^3\oplus\framebox{${\mathbb{R}}^3$}
&&{\mathbb{R}}^3
\end{array}\end{equation}
where ${\mathrm{S}}^2{\mathbb{R}}^2$ denotes symmetric $3$-tensors whilst skew 
$3$-tensors are identified with ${\mathbb{R}}^3$ using $\epsilon_i{}^{jk}$. In
this diagram, the spaces indicated thus
\raisebox{2pt}{\framebox{\raisebox{5pt}{{\quad}}}} are joined by isomorphisms
indicated thus~$\nearrow$. A simple diagram chase cancels these spaces and
results in the linear elasticity complex (\ref{linearelasticitycomplex}).
In~\cite{bulletin}, Arnold, Falk, and Winther use a halfway-house complex
$$\begin{array}{ccccccc} {\mathbb{R}}^3&&&&&&{\mathbb{R}}^3\\
\oplus&\longrightarrow&
{\mathbb{R}}^3\oplus{\mathrm{S}}^2{\mathbb{R}}^3&\longrightarrow&
{\mathrm{S}}^2{\mathbb{R}}^3\oplus{\mathbb{R}}^3&\longrightarrow&\oplus\\
{\mathbb{R}}^3&&&&&&{\mathbb{R}}^3
\end{array}$$
obtained by cancelling only
$\framebox{${\mathbb{R}}^3\otimes{\mathbb{R}}^3$}\nearrow
\framebox{${\mathbb{R}}^3\otimes{\mathbb{R}}^3$}$ from~(\ref{bigdiagram}), to
construct new and stable finite element schemes for linear elasticity mimicking
the previously known stable finite element schemes for the grad-curl-div 
complex.

\subsection*{Projective geometry}
The connection (\ref{tractors}) may be viewed as follows. Consider the unit
three-sphere $S^3\subset{\mathbb{R}}^4$. There is no difficulty in taking the
gradient of a function $F$ on $S^3$ with values in the skew
$2$-tensors~$\Lambda^2{\mathbb{R}}^4$. However, each point on $S^3$ is also a 
vector $v\in{\mathbb{R}}^4$, which may be used to decompose these skew 
$2$-tensors:--
$$\Lambda^2{\mathbb{R}}^4=
\{\omega\mbox{ s.t.\ }v\wedge\omega=0\}\oplus
\{\omega\mbox{ s.t.\ }v\intprod\omega=0\}
\cong{\mathbb{R}}^3\oplus{\mathbb{R}}^3.$$
This decomposition does not see the sign of $v$ and so descends to the quotient
of $S^3$ under antipodal identification, namely real projective
$3$-space~${\mathbb{RP}}_3$. The upshot is that the gradient of $F$ may be
written in terms of the intrinsic calculus on ${\mathbb{RP}}_3$ and viewed in a
standard affine co\"ordinate patch
${\mathbb{R}}^3\hookrightarrow{\mathbb{RP}}_3$. The result is~(\ref{tractors}).
The construction of (\ref{linearelasticitycomplex}) from (\ref{bigdiagram}) is
due to Calabi~\cite{calabi}. It may also be viewed as a geometric realisation
of the Jantzen-Zuckerman translation principle from representation
theory~\cite{vogan} and, as such, admits vast generalisation in the newly
developed field of parabolic geometry~\cite{thebook}.

\renewcommand{\section}{\subsection}

\end{document}